%% file: paper.tex
\documentclass[a4paper,abstracton]{scrartcl}
\usepackage{graphicx}
\usepackage{amsmath}
\usepackage{amsfonts}
\usepackage{booktabs}
\usepackage{comment}
\usepackage{tablefootnote}
\usepackage{xcolor}
\usepackage{tikz}

\usepackage{authblk}

\begin{document}

\title{A new formulation for the collection and delivery problem of biomedical specimen}

\author[1]{Luis Aurelio Rocha}
\author[1]{Alena Otto}
\author[2]{Marc Goerigk}

\affil[1]{Chair of Management Science / Operations and Supply Chain Management, University of Passau, Germany}
\affil[2]{Chair of Business Decisions and Data Science, University of Passau, Germany}

\date{}

\maketitle

\begin{abstract}
We study \textit{the \underline{c}ollection and \underline{d}elivery problem of biomedical \underline{s}pecimens (CDSP)} with multiple trips, time windows, a homogeneous fleet, and the objective of minimizing total completion time of delivery requests. This is a prominent problem in healthcare logistics, where  specimens (blood, plasma, urin etc.) collected from patients in doctor's offices and hospitals are transported to a central laboratory for advanced analysis.
To the best of our knowledge, available exact solution approaches for CDSP have been  able to solve only small instances with up to 9 delivery requests. In this paper, we propose a two-index mixed-integer programming formulation that, when used with an off-the-shelf solver, results in a fast exact solution approach.  Computational experiments on a benchmark data set confirm that the proposed formulation outperforms both the state-of-the-art model and the state-of-the-art metaheuristic from the literature, solving 80 out of 168 benchmark instances to optimality, including a significant number of instances with 100 delivery requests.
\end{abstract}

\noindent\textsf{\textbf{Keywords:}} Biomedical specimen; Routing; Multitrip; Healthcare logistics; MIP

\noindent\textsf{\textbf{Acknowledgments:}} This study was funded by the German Federal Ministry for Digital and Transport (BMDV) as part of project KIMoNo (45KI01A011).

\section{Introduction}
\label{sec.Intro}

Most medical decisions rely on specialized laboratory tests, which require a relevant specimen or sample (e.g., blood, plasma, or urine) to be taken from the patient and sent to a medical laboratory. 
In \textit{regular} collection of biomedical specimens, vehicles visit various points of care, such as doctor's offices or clinics, once a day to collect specimens (which form a \textit{delivery request}) and deliver them to the laboratory.

Specimen collection and delivery differs from typical vehicle routing problems  in several aspects. First, capacity constraints are negligible  since specimens take up little space. Second, variable transportation costs (e.g., fuel) are less significant   and usually not factored into the problem. Third, as many specimens as possible should arrive at the laboratory as early as possible to ensure an earlier finalization of specimen  analytics. Indeed, the laboratory capacity is limited and processing in the laboratory can only begin \textit{after} specimen arrival. Instead of bringing all specimens at once, as with a makespan objective, it is more suitable to minimize the sum of request completion times. Here, the request completion time is the time spent until the arrival of the delivery request at the laboratory. This approach ensures many specimens arrive earlier, even if it delays the delivery of a few others. As a result, the laboratory can start processing specimens earlier, allowing for the analysis of all specimens to be completed sooner.

\textit{The \underline{c}ollection and \underline{d}elivery problem of biomedical \underline{s}pecimen (CDSP)} with the objective of minimizing total completion time of delivery
requests was first proposed by Zabinsky et al. \cite{Zabinsky.2020},
who suggested a mixed-integer programming formulation and a customized branch-and-bound algorithm \textit{VeRSA}, which can solve instances with up to 9 requests. In a recent study, Ferone et al. \cite{Ferone.2023} designed a problem-specific \textit{\underline{a}daptive \underline{l}arge \underline{n}eighborhood \underline{s}earch (ALNS)} tested against VeRSA on instances with 9, 25, and 57 requests. VeRSA solved instances with 9 requests within a one-hour time limit per instance and found feasible solutions only in 60\% of cases, all feasible solutions found by VeRSA for larger instances with 57 requests had bad objective values (at least 12 times bigger than that found by ALNS). ALNS was able to find all optimal solutions reached by VeRSA and significantly improved solutions in all other cases. The model and algorithms developed by \cite{Zabinsky.2020} and \cite{Ferone.2023} remain the state-of-the-art approaches to CDSP to date.

Boland et al. \cite{Boland.2000} introduced a replenishment-arc model formulation for the constrained traveling salesman problem, which has since been successfully applied to multi-trip routing problems (\cite{Neira.2020,Rivera.2016}, see Cattaruzza et al. \cite{Cattaruzza.2016} for a review).

In this paper we demonstrate that the replenishment-arc formulation can be adapted to CDSP, resulting in the best-known exact solution approach for this problem so far. We outline the problem and the suggested two-index model formulation in Section \ref{sec.ProblemDef}.  Section \ref{sec.CompExp} presents computational experiments, and we conclude with a discussion and future research directions in Section \ref{sec.Discussion}.

\section{Problem definition and new model formulation}
We formulate CDSP as presented in \cite{Zabinsky.2020}. 
Directed multigraph $G(V,A,c)$ has the set of nodes $V=P\cup\{0\}$ consisting of points of care $P=\{1,\ldots,n\}$ and \textit{depot} 0, which denotes the laboratory. Each point of care generates a specimen collection \textit{request}, which needs to be transported to the depot. In the following, we refer to set $P$  as the set of points of care and the set of requests interchangeably, depending on the context. The set of arcs $A=A^0\cup A^P\cup A^R$ consists of  \textit{original} arc sets $A^0$ and $A^P$ and of an artificial \textit{replenishment arc set} $A^R$. Arc labels $c_e, \forall e\in A,$ denote the  travel time to traverse arc $e$. We assume that the subgraph $G'=(V,A^0\cup A^P, c)$ spanned over the original arcs is a complete directed graph with no multiple edges and no loops, for which triangle inequalities hold: $c_{ij}+c_{jk}\geq c_{ik}, ~\forall (i,j), (j,k), (i,k)\in A^0\cup A^P$. Arcs $e\in A^0$ denote the arcs \textit{adjacent to the depot} and set $A^P$ refers to the remaining original arcs.

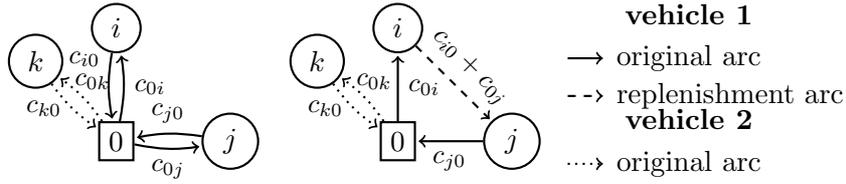
\begin{figure}[bt]
\centering
\input{example_replenishment_arc_schema}
\caption{Vehicle tours using only original arcs (left) and using replenishment arcs (right)}
\label{fig.replenishment_arcs}
\end{figure}

The $K\in \mathbb{N}$ vehicles start and end their tours at the depot. A part of a vehicle's tour between two consecutive depot visits is called a \textit{trip}. CDSP is a {multi-trip} routing problem, as each vehicle may visit the depot several times to deliver collected specimens. 
Replenishment arcs $A^R$ are a key element of the suggested problem formulation. They  transform CDSP into a \textit{single-trip} routing problem by 'hiding' depot visits within replenishment arc, as shown in Figure~\ref{fig.replenishment_arcs}. We introduce exactly one replenishment arc between each ordered pair of distinct nodes $i,j\in P,i\neq j.$ Let $c_{i0}$ and $c_{0j}$ be the travel times from point of care $i$ to the depot and from the depot to point of care $j, j\neq i,$ using the respective original arcs, then the arc traversal time of replenishment arc $e$ spanned between $i$ and $j$ is set to $c_e=c_{i0}+c_{0j}$.

\label{sec.ProblemDef}

The \textit{(shift)} of each vehicle, or the time between leaving the depot for the first trip and returning to the depot after the last trip, cannot exceed $\tau_{max}$. There is a hard time window $[r_j, d_j]$ associated with each point of care $j\in P$. If a vehicle arrives at $j$ before release time $r_j$, it has to wait. Moreover, the deadline $d_0$ restricts the return time of the last vehicle to the depot. 

The \textit{completion time} of request $i\in P$ denotes its delivery time to the laboratory.
The \textit{objective} of CDSP is to find vehicle tours in $G$ that start and end at the depot, minimizing the sum of all  request completion times. Thereby, \textit{(i)} each point of care $i\in P$ is visited exactly once within its time window; \textit{(ii)}  at most $K$ vehicles leave the depot; \textit{(iii)} 
each vehicle's shift does not exceed $\tau_{max}$.

We define the following decision variables. Binary variable $x_{e}\in\{0,1\}, \forall e\in A,$ equals 1 if some vehicle traverses arc $e$, and 0 otherwise. Nonnegative variable $z_{i}\in \mathbb{R}^+_0$ tracks the time each point of care $i\in P$ is visited by a vehicle.
Variable $C_j\in \mathbb{R}^+_0, \forall j\in P,$ corresponds to the completion time of  request $j$. To compute $C_j$,  we define binary variables $y_{i,j}\in\{0,1\}, \forall i,j\in P$, which equal 1 if node $i$ is visited by the same vehicle after collecting request $j$ and before returning to the depot. Since points of care cannot be visited before their release times, some vehicles may start their shift with a delay. Therefore, to compute the shift duration, we introduce variables $\tau_{j}\in \mathbb{R}, \forall j\in P,$ where $\tau_{j}$ denotes the time elapsed from the beginning of the vehicle's shift until it visits $\tau_{j}$.   

To simplify notation, we denote the source and the target of each arc $e\in A$ as $s(e)$ and $t(e)$, correspondingly. We further define
$\delta^+_j:=\{e\in A | s(e)=j\}$
and
$\delta^-_j:=\{e\in A| t(e)=j\}$
 as the sets of outgoing and incoming arcs of node $j\in V$, correspondingly. In the \textit{preprocessing}, we also update the release times and deadlines as $r_j:=\max\{r_j,c_{0j}\}$ and $d_j:=\min\{d_j,d_0-c_{j0}\}$, $\forall j\in P$.

\begingroup
\hspace{1.0cm}
\allowdisplaybreaks
\begin{align}
\label{eq.SCP.CompletionTime.obj}
\text{\textit{Minimize }} & F(x, y,\tau, z,C) = \sum_{j \in P} C_j \\
\label{eq.SCP.startEndNodeA}
\text{ s.t.} \qquad \qquad  \sum_{e \in \delta^+_0} x_{e} &= \sum_{e' \in \delta^-_0} x_{e'} \leq K \\
\label{eq.SCP.CustomerJustOnce}
\sum_{e \in \delta^+_j} {x_{e}} &=  \sum_{e' \in \delta^-_j} {x_{e'}}=1  \quad \quad \forall j \in P\\
\label{eq.SCP.timeFlow}
  z_{s(e)} + c_{e} & \leq z_{t(e)} +M_{e}(1-x_{e}) \quad \quad \forall e \in A^P\cup A^R\\
\label{eq.SCP.HardTimeWindows}
r_j \leq z_{j}  &\leq d_{j}  \quad \quad \forall j \in P\\
\label{eq.SCP.requestCollectionStart}
y_{jj}& =1  \ \forall j \in P\\
\label{eq.SCP.requestContinuationA}
y_{s(e),j} \leq y_{t(e),j}  +(1 &- x_{e})  \: \forall e \in A^P, j \in P, j\neq t(e)\\
\label{eq.SCP.CompletionTimeConstraint}
z_{i}+ c_{i0}  & \leq C_j +M_i(1-y_{ij}) \quad \quad \forall i,j \in P\\
 \label{eq.SCP.tripDurationStartNewTripB}
 c_{0j} & \leq \tau_{j}  \quad \quad \forall j \in P\\
\label{eq.SCP.tripDurationContinuation}
\tau_{s(e)}+(z_{t(e)} -z_{s(e)})  & \leq \tau_{t(e)} +M'_e(1-x_{e}) \quad \quad \forall e \in A^P\cup A^R\\
\label{eq.SCP.maxTripDuration}
\tau_{j}+c_{j0} & \leq \tau_{max}  \quad \quad \forall j \in P\\
\label{eq.SCP.BinaryDomainsA}
x_{e} & \in \{0,1\} \quad \quad \forall e\in A\\
\label{eq.SCP.BinaryDomainsC}
y_{ij} & \in \{0,1\} \quad \quad \forall i,j\in P\\
\label{eq.SCP.ContDomains}
z_{j}, \tau_{j}, C_j & \geq 0 \quad \quad \forall j\in P
\end{align}
\endgroup

Objective function \eqref{eq.SCP.CompletionTime.obj} minimizes the total request completion time. 
Constraints~\eqref{eq.SCP.startEndNodeA} and \eqref{eq.SCP.CustomerJustOnce} are flow conservation constraints. At most $K$ vehicles leave the depot and each point of care $j\in P$ is visited exactly once. Constraints~\eqref{eq.SCP.timeFlow} compute the visit times of points of care, thereby $M_{e}:= \max\{0,d_{s(e)}+c_{e}-r_{t(e)} \} ~\forall e\in A^P\cup A^R$. Constraints~\eqref{eq.SCP.HardTimeWindows} ensure that vehicles visit each point of care within the respective time window. Request $j$ is collected, when point of care $j$ is visited (Constraints~\eqref{eq.SCP.requestCollectionStart}). Constraints~\eqref{eq.SCP.requestContinuationA} mark the nodes visited by request $j$ before it reaches the depot and Constraints~\eqref{eq.SCP.CompletionTimeConstraint} compute completion times of each request. We set $M_i:=d_i+c_{i0} ~\forall i\in P$. Constraints~\eqref{eq.SCP.maxTripDuration} restrict the  shift duration of each vehicle, and the respective variables $\tau_j$ are initialized in Constraints~\eqref{eq.SCP.tripDurationStartNewTripB}. Constraints~ \eqref{eq.SCP.tripDurationContinuation} compute the visit time of each point of care relative to the shift start of the respective vehicle. We set $M'_e:=d_{t(e)}-c_{0,t(e)} ~\forall e\in A^P\cup A^R.$ Finally, Constraints~\eqref{eq.SCP.BinaryDomainsA}-\eqref{eq.SCP.ContDomains} describe the domains of the variables.

\section{Computational experiments}
\label{sec.CompExp}

We conducted computational experiments on the benchmark data set of \cite{Ferone.2022}, which is based on the well-known Solomon instances. The dataset consists from 168 instances grouped in $3\times 3\times 2=18$ settings which differ in  
\begin{itemize}
\item the instance size: small ($n=25$), medium ($n=50$), and large ($n=100$) instances have 10, 15, and 25 vehicles, respectively;
\item the spatial distribution of nodes: \textit{clustered [C]}, uniformly at \textit{random} over a square \textit{[R]}, and a \textit{combination} of random and clustered \textit{[RC]};
\item the width of the \textit{time windows (TW)}: wide and tight.
\end{itemize}
Following Zabinsky et al. \cite{Zabinsky.2020} and Ferone et al. \cite{Ferone.2023}, who align their reporting to a common key performance indicator, we subtract a constant ($\sum_{j\in P}{r_j}$) from the original objective value and report $F'=\sum_{j\in P}{(C_j-r_j)}$.
 
 We implemented the model of Zabinsky et al. \cite{Zabinsky.2020} (called simply as \textit{Model \cite{Zabinsky.2020}} in the following) and our model on the same machine using Gurobi Optimizer 11.0.3 with standard parameters, limited to 16 threads, and a time limit of 3600 seconds per model and instance. The experiments were conducted on a Linux server (6248R CPU, 3.00GHz$\times$96 processors, 754.5 GB RAM). Since Model \cite{Zabinsky.2020} could not find a feasible solution for any small instance ($n=25$) within the time limit, we did not run it on larger instances. 

Table~\ref{tab:results} compares the results of our model, Model \cite{Zabinsky.2020} and ALNS \cite{Ferone.2023}. The last two columns report average objective value and the total time to the best found solution \textit{(Ttb)} as reported in \cite{Ferone.2023}.  For $n=25$,  Ferone et al. \cite{Ferone.2023} report only the results over 26 and 18 instances out of 56 for Avg $F'$ and Ttb[s], respectively. We mark these cases with * in the table.
Our model clearly outperformed Model~\cite{Zabinsky.2020}. It solved 55 out of 56 small instances ($n=25$), whereas Model~\cite{Zabinsky.2020} was not able to find a feasible solution to any of these instances. 

For medium ($n=50$) and large ($n=100$) instances, our model was able to find and prove optimal solutions in 30\% (17 out of 56) and 14\% (8 out of 56) of cases, respectively. It outperformed ALNS~\cite{Ferone.2023} in the average achieved objective values in eight out of 12 settings. The gaps achieved by our model remain moderate at 1.12\% on average and never exceed 14\%.

\begin{table}[htbp]
\setlength{\tabcolsep}{0.3em}
\centering
\caption{Comparative performance of models and algorithms}
\label{tab:results}
\begin{tabular}{llllllllllll}
\hline
& & & \multicolumn{5}{c}{Our model}&  \multicolumn{2}{c}{Model \cite{Zabinsky.2020}}  & \multicolumn{2}{c}{ALNS \cite{Ferone.2023}} \\
\cmidrule(lr){4-8} \cmidrule(lr){9-10} \cmidrule(lr){11-12}
\multicolumn{3}{l}{Instance} & Avg & Avg & Avg  & \# & \# & Avg & Avg & Avg & Avg \\
$n$ & class & TW & $F'$  & gap & T[s] &opt & TL & $F'$ & T[s] & $F'$ & Ttb[s] \\
\hline
25 & C & tight & \textbf{672} & 0.00~\% & 168&9 & - & - & 3600 & 755* & 482* \\ 
 & & wide & \textbf{789} & 0.00~\% & 21&8 & - &  - & 3600 & 801* & 357* \\ 
 & R & tight & \textbf{829} & 0.00~\% & 42&12 & - & - & 3600 & 856* & 250* \\ 
 & & wide & \textbf{835} & 0.00~\%& 464 &11& - & - & 3600 & 802* & 645*\\ 
 & RC & tight & \textbf{1193} & 0.00~\% & 53&8 & - & - & 3600 & 1359* & 527* \\ 
 & & wide & \textbf{1134} & 0.01~\%& 864&7& 1& - & 3600 & 1148* & 820*\\ 
 50 & C & tight & \textbf{1467} & 0.34~\%& 2403& 3 & 6 &   -&-&1468 & 1290 \\ 
 & & wide & \textbf{1708}& 0.17~\%& 1354& 5 & 3 &   -&-&1709& 1362 \\ 
 & R & tight & \textbf{1897}& 1.91~\%& 2178& 5 & 7 &   -&-&1943 & 1147 \\ 
 & & wide & \textbf{1870}& 1.26~\%& 2947& 2 & 9 &   -&-&1872 & 1374 \\ 
 & RC & tight & \textbf{2737}& 1.46~\%& 3234& 1 & 7 &   -&-&2743 & 1345 \\ 
 & & wide & \textbf{2528} & 0.46~\%& 3152& 1 & 7 &   -&-&2531 & 1526 \\ 
100 & C & tight & 3708& 1.71~\%& 3245& 1 & 8&   -&-&\textbf{3545}& 2064 \\ 
 & & wide & 4389& 1.43~\%& 2801& 2 & 6&   -&-&\textbf{3970}& 2440 \\ 
 & R & tight & \textbf{3552}& 2.95~\%& 3333& 1 & 11 &   -&-&3693 & 1732 \\ 
 & & wide & 3745& 3.13~\%& 3011& 2 & 9&   -&-&\textbf{3529} & 2289 \\ 
 & RC & tight & \textbf{4506}& 2.93~\%& 3601& - & 8 &   -&-&4694 & 2178 \\ 
 & & wide & 4490& 2.36~\%& 2848& 2 & 6&   -&-&\textbf{4223}& 2498 \\ 
\hline
\multicolumn{12}{l}{\parbox{10cm}{\footnotesize{$F'=\sum_{j\in P}{(C_j-r_j)}$ is the total completion time corrected by a constant (release time),}}}\\ 
\multicolumn{12}{l}{\footnotesize{T[s] -- runtime in seconds,}}\\ 
\multicolumn{12}{l}{\footnotesize{\# opt -- \# instances with optimal solutions found and proven, }}\\
\multicolumn{12}{l}{\footnotesize{\# TL -- \# instances with reached time limit,}}\\
\multicolumn{12}{l}{\footnotesize{Ttb[s] -- time to the best solution in seconds,}}\\
\multicolumn{12}{l}{\footnotesize{The best average objective values are marked in \textbf{bold}. }}\\
\multicolumn{12}{l}{\footnotesize{Partial averages are marked with *.}}
\end{tabular}
\end{table}

\section{Discussion and future research directions}
\label{sec.Discussion}

The presented two-index model formulation outperforms the state-of-the-art approaches to CDSP.
Further research may develop further efficient customized heuristic and exact solution approaches for CDSP. The problem can be extended to multiple laboratories and include fleet dimensioning by introducing fixed vehicle costs to the objective.

\end{document}

%% file: example_replenishment_arc_schema.tex
\begin{tikzpicture}[->,shorten >=1pt,auto,node distance=1.5cm,
                    thick,cust node/.style={circle,draw},delivery node/.style={rectangle,draw},dummy node/.style={rectangle,fill=white}]
 
 \node[delivery node] (4) {$0$};
  \node[cust node] (5) [above of=4] {$i$};
  \node[cust node] (6) [right of=4] {$j$};
  \node[cust node] (7) [above left of=4] {$k$};
 \node[delivery node] (0) [right of=6,xshift=20pt] {$0$};
  \node[cust node] (1) [above of=0] {$i$};
  \node[cust node] (2) [right of=0] {$j$};
  \node[cust node] (3) [above left of=0] {$k$};

  \path[every node/.style={font=\sffamily\small}]
    (0) edge node[right]  {$c_{0i}$} (1)
    (0) edge[dotted, bend right=10] node[right,pos=0.9]  {$c_{0k}$} (3)
    (3) edge[dotted, bend right=10] node[left]  {$c_{k0}$} (0)
    (1) edge[dashed] node[above,sloped] {$c_{i0}+c_{0j}$} (2)
    (2) edge node[below] {$c_{j0}$} (0)
    (4) edge[bend right=10] node[right]  {$c_{0i}$} (5)
    (5) edge[bend right=10] node[left,pos=0.08]  {$c_{i0}$} (4)
    (4) edge[dotted, bend right=10] node[right,pos=0.9]  {$c_{0k}$} (7)
    (7) edge[dotted, bend right=10] node[left]  {$c_{k0}$} (4)
    (6) edge[bend right=10] node[above] {$c_{j0}$} (4)
    (4) edge[bend right=10] node[below] {$c_{0j}$} (6);

\draw (1.east) ++(1.9,-0.4) coordinate (tmp) -- +(0.5,0) node[right,label=above:{\textbf{vehicle 1}}] {original arc};
\draw [dashed] (tmp)  ++(0,-0.5) -- +(0.5,0) node[right] {replenishment arc};
\draw [dotted] (tmp)  ++(0,-1.4) -- +(0.5,0) node[right,label=above:{\textbf{vehicle 2}}] {original arc};
\end{tikzpicture}